\documentclass[manuscript, printscheme]{aam}

\usepackage{stmaryrd}

\DeclareMathOperator*{\Loss}{Loss}

\def\bx{\mathbf{x}}
\def\bX{\mathbf{X}}

\def\bn{\mathbf{n}}
\def\bftau{{\boldsymbol\tau}}
\def\bb{\mathbf{b}}
\def\bu{\mathbf{u}}
\def\bv{\mathbf{v}}
\def\bff{\mathbf{f}}
\def\bg{\mathbf{g}}
\def\bU{\mathbf{U}}
\def\bF{\mathbf{F}}
\def\bI{\mathbf{I}}
\def\bft{{\boldsymbol\theta}}
\def\grad{\nabla}
\def\pd{\partial}
\newcommand{\beq}{\begin{equation}}
\newcommand{\eeq}{\end{equation}}
\newcommand{\beqs}{\begin{eqnarray}}
\newcommand{\eeqs}{\end{eqnarray}}
\newcommand{\beqsn}{\begin{eqnarray*}}
\newcommand{\eeqsn}{\end{eqnarray*}}
\newcommand{\bary}{\begin{array}}
\newcommand{\eary}{\end{array}}

\def\div{\nabla\cdot}
\newcommand*{\abs}[1]{\left|#1\right|}
\newcommand*{\dbblk}[1]{{\bigl\llbracket} {#1} {\bigr\rrbracket}}
\newcommand*{\dbblkL}[1]{{\biggl\llbracket} {#1} {\biggr\rrbracket}}

\begin{document}
\title{A discontinuity and cusp capturing PINN for Stokes interface problems with discontinuous viscosity and singular forces}


\author[Y.-H. Tseng et~al.]{Yu-Hau Tseng\affil{1} and Ming-Chih Lai\affil{2,}\corrauth}
\address{\affilnum{1}\ Department of Applied Mathematics, National University of Kaohsiung, Kaohsiung 81148, Taiwan. \\ \affilnum{2}\ Department of Applied Mathematics, National Yang Ming Chiao Tung University, 1001 Ta Hsueh Road, Hsinchu 300, Taiwan.}
 \emails{{\tt mclai@math.nctu.edu.tw} (M.-C. Lai)}

\begin{abstract}
In this paper, we present a discontinuity and cusp capturing physics-informed neural network (PINN) to solve Stokes equations with a piecewise-constant viscosity and singular force along an interface. We first reformulate the governing equations in each fluid domain separately and replace the singular force effect with the traction balance equation between solutions in two sides along the interface. Since the pressure is discontinuous and the velocity has discontinuous derivatives across the interface, we hereby use a network consisting of two fully-connected sub-networks that approximate the pressure and velocity, respectively. The two sub-networks share the same primary coordinate input arguments but with different augmented feature inputs. These two augmented inputs provide the interface information, so we assume that a level set function is given and its zero level set indicates the position of the interface. The pressure sub-network uses an indicator function as an augmented input to capture the function discontinuity, while the velocity sub-network uses a cusp-enforced level set function to capture the derivative discontinuities via the traction balance equation. We perform a series of numerical experiments to solve two- and three-dimensional Stokes interface problems and perform an accuracy comparison with the augmented immersed interface methods in literature. Our results indicate that even a shallow network with a moderate number of neurons and sufficient training data points can achieve prediction accuracy comparable to that of immersed interface methods.
\end{abstract}


\maketitle

\section{Introduction}
\label{sec:intro}
The incompressible two-phase flow problems have a wide range of applications in various scientific and engineering fields; see \cite{JR93} and the references therein. Very often, this kind of problem involves solving Stokes equations with a discontinuous viscosity and singular force (such as surface tension force) along the fluid interface. For simplicity, we call this problem the Stokes interface problem hereafter. Solving Stokes interface problems numerically is known to be quite challenging in the scientific computing community, giving rise to several difficulties. One major difficulty is to handle the coupling between two adjacent fluids and the interfacial boundary conditions (arising from the singular force). Another difficulty comes from that the pressure is, in fact, discontinuous, and the partial derivatives of the velocity are also discontinuous across the interface. Thus, careful numerical treatments must be adopted for accurate discretizations near the interface, regardless of using the finite difference or finite element method.

In order to address the above issues, Peskin proposed the immersed boundary (IB) method~\cite{Peskin02}, which introduces a regularized version of the discrete delta function to discretize the singular force in a Cartesian grid setting. The discontinuous viscosity can also be regularized by a smoothed version of the Heaviside function. In this framework, the problem becomes smooth so that regular finite difference discretizations such as the MAC scheme can be applied to solve the Stokes equations. However, it is known that the IB method achieves only first-order accuracy~\cite{LM12} for the solution variables. To improve the accuracy, an alternative approach is to solve the Stokes equations within each fluid domain separately and reformulate the singular force term into the traction balance equation on the interface. The above two formulations will be given in Section 2. In this paper, we shall introduce a neural network methodology to solve the problem based on the second formulation. Since the present goal here is to make comparisons with traditional methods such as immersed interface method (IIM, also based on the second formulation) from the implementation and accuracy aspects, we provide a brief review of the IIM in the following.

The immersed interface method proposed by LeVeque and Li~\cite{LL94} was originally designed to solve elliptic interface problems with discontinuous coefficients. Their approach incorporates the derived jump conditions via local coordinates into the finite difference discretization so that the local truncation error can be reduced appropriately to achieve the solution in second-order accuracy in the $L^\infty$ norm. Lai et al.~\cite{HLY15,LT08} later proposed a simplified version of the IIM that directly utilizes jump conditions to achieve second-order accuracy in the $L^\infty$ norm without resorting to local coordinates. Of course, there are other variants of solving elliptic interface problems using the jump conditions, such as the ghost fluid method \cite{EG20,LFK00}, or meshless method \cite{Oru21}. However, it is not our intention to give an exhaustive review here so readers who are interested in the recent results can refer to \cite{Oru21}, and the references therein. Comparing with the well-developed solvers for elliptic interface problems, there are relatively few works to investigate the present Stokes interface problems with discontinuous viscosity and singular forces. As we mentioned earlier, the major difficulty is that the solution and derivative discontinuities are all coupled, making it hard to develop accurate numerical methods. In \cite{LIL07}, Li et al.~developed a new IIM-based method for the two-dimensional problem that decouples the jump conditions of the fluid variables by introducing two augmented velocity variables. Although the new augmented system of equations can be respectively written as three Poisson interface problems for the pressure and the augmented velocity, they are still coupled through the augmented variable jumps defined on the interface. As a result, the GMRES iterative method is used to solve the Schur complement system for the augmented variable jumps. Recently, Wang et al.~\cite{WT22} used the similar idea but a simple version of IIM and extended the method to solve the three-dimensional problems.

In addition to traditional grid-based methods, the scientific computing community has shown increasing interest in utilizing deep neural networks to tackle elliptic interface problems. One advantage of the neural network approach is its mesh-free nature, enabling us to handle problems with complex interfaces or irregular domains more easily. However, with the common usage of the smooth activation function in neural network methods, the network solution is naturally smooth and cannot capture the solution behaviors (function discontinuity or derivative discontinuity) sharply. One simple way to circumvent  such difficulty is to decompose the computational domain into two sub-domains and employs separate networks for each sub-domain such as the networks in \cite{GY22,HHM22,WL22}. Recently, the authors have developed a discontinuity capturing shallow neural network to solve the elliptic interface problems \cite{HLL22}. The main idea is to represent a $d$-dimensional discontinuous function by a single ($d+1$)-dimensional continuous neural network function in which the extra augmented variable labels the sub-domains. Once the function representation is chosen, the rest of the method follows the physics-informed neural network (PINN) framework developed in \cite{RPK19}. Meanwhile, we have also developed a ($d+1$)-dimensional cusp capturing neural network \cite{TLHL22} to represent a $d$-dimensional continuous function but with discontinuous partial derivatives across the interface. In such a case, the extra augmented variable uses the cusp-enforced level set function instead. As we mentioned earlier, in the present Stokes interface problem, the pressure is usually discontinuous, and the partial derivatives of velocity are discontinuous, so they are well represented by the discontinuity capturing and cusp capturing networks, respectively. Both network function representations will be given in detail in Section~\ref{sec:CuspDCNN}.

The remainder of the paper is organized as follows. In Section~\ref{sec:StokesInterfaceEqns}, we first present the Stokes equations with a discontinuity viscosity and singular force. Then we reformulate the governing equations in each sub-domain and introduce the traction balance equation to replace the singular force effect. In Section~\ref{sec:CuspDCNN}, we introduce the discontinuity and cusp capturing PINN for solving the model problems. Numerical experiments demonstrating the accuracy of the proposed method compared with the existing augmented immersed interface methods are shown in Section~\ref{sec:numerics}. Some concluding remarks are given in Section~\ref{sec:conclusion}.

\section{Stokes equations with a discontinuous viscosity and singular force}
\label{sec:StokesInterfaceEqns}
In this paper, we aim to solve a $d$-dimensional incompressible Stokes system with a piecewise-constant viscosity and singular force on an embedded interface. This problem arises from the simulation of the flow of two adjacent fluids with different viscosity in the presence of interfacial force. We start by writing down the governing equations in the immersed boundary (IB) formulation \cite{Peskin02}, which regards the interface as a singular force generator so that the equations can be described in the whole fluid domain $\Omega$. To begin with, let $\Omega\subset\mathbb{R}^d$, where $d=2$ or $3$, be a bounded domain, and $\Gamma$ be an embedded co-dimensional one $C^1$-interface that separates $\Omega$ into two subdomains, $\Omega^-$ and $\Omega^+$, such that $\Omega = \Omega^- \cup \Omega^+ \cup \Gamma$. This Stokes interface problem can be written as in \cite{LIL07,WT22}
	\beqs
	\nabla p &=& \nabla\cdot\left(\mu\left(\nabla{\bu}+\nabla{\bu}^T\right)\right) + \bff + {\bg}, \quad \mbox{in} \, \Omega, \label{eq:Stokes_deltaform}\\
	\nabla\cdot{\bu} & = & 0, \quad \mbox{in} \, \Omega, \label{eq:div0}\\
	{\bu}  & = & {\bu}_b, \quad \mbox{on} \, \partial\Omega. \label{eq:BCs_D0}
	\eeqs
Here, $\bu=(u_1,\ldots,u_d)\in\mathbb{R}^d$ is the fluid velocity, $p$ is the pressure, and $\bg=(g_1,\ldots,g_d)$ is an external force, which are all functions of $\bx=(x_1, \ldots x_d) \in \Omega$. The viscosity $\mu(\bx)$ is defined piecewise-constantly such that $\mu(\bx) = \mu^\pm$ if $\bx \in \Omega^\pm$. The force term $\bff$ is a singular concentrated body force expressed in the IB formulation as $\bff(\bx) =\int_\Gamma \bF(\bX)\delta(\bx-\bX)dS$, where the interfacial force $\bF=(F_1,\ldots,F_d)$ and $\delta(\bx)$ is the $d$-dimensional Dirac delta function. Due to the delta function singularity, one can immediately see that the pressure and velocity should not be smooth mathematically. In fact, one can reformulate the governing equations (\ref{eq:Stokes_deltaform})-(\ref{eq:BCs_D0}) in each fluid domain  $\Omega^\pm$ separately, but the solutions are linked by the traction balance equation (\ref{eq:jump_stress0}) on the interface $\Gamma$ as follows.
	\beqs
		-\nabla p + \mu\Delta{\bu} + {\bg}&=& {\bf 0}, \quad \mbox{in} \quad\Omega^\pm, \label{eq:Stokes_jumpform}\\
		\div{\bu} & = & 0, \quad \mbox{in} \quad\Omega^\pm, \label{eq:div1} \\
		\dbblk{-p\bI + \mu\left(\nabla{\bu}+\nabla{\bu}^T\right)}\cdot\bn+{\bF} &=& \mathbf{0}, \quad \mbox{on} \quad \Gamma, \label{eq:jump_stress0} \\
		\bu &=& \bu_b, \quad \mbox{on} \quad \partial\Omega.\label{eq:BCs_D}	
	\eeqs
Here, the notation $\dbblk{\cdot}$ represents the jump of a quantity across the interface; that is,  for any $\bx_\Gamma\in\Gamma$,
	\beq
	\dbblk{f}(\bx_\Gamma) = \lim_{\bx\in\Omega^+,\,\bx\rightarrow\bx_\Gamma} f(\bx) - \lim_{\bx\in\Omega^-,\,\bx\rightarrow\bx_\Gamma}f(\bx) = f^+(\bx_\Gamma) - f^-(\bx_\Gamma). \label{eq:jump_def}
	\eeq
The unit outward normal vector, $\bn=(n_1,\ldots,n_d)$, points from $\Omega^-$ to $\Omega^+$ along the interface $\Gamma$. This above traction  balance equation (\ref{eq:jump_stress0}) can be derived  in \cite{Poz97}.

Since the fluid is viscous, the velocity is continuous across the interface. However, the traction balance equation (\ref{eq:jump_stress0}) shows that the pressure is discontinuous and the velocity has discontinuous partial derivatives. Indeed, we can rewrite Eq.~(\ref{eq:jump_stress0}) component-wisely as
	\beqs
	 -\dbblk{p}n_k + \dbblk{\mu\frac{\pd u_k}{\pd \bn}} + \dbblkL{\mu\frac{\pd\bu}{\pd{x}_k}}\cdot\bn +F_k=0, \quad k=1,\ldots,d. \label{eq:jump_stressfxyz}
	\eeqs
In the next section, we shall present a neural network learning methodology based on PINN \cite{RPK19} to solve Eqs.~(\ref{eq:Stokes_jumpform})-(\ref{eq:div1}) with the traction balance equation (\ref{eq:jump_stressfxyz}), subject to the boundary condition (\ref{eq:BCs_D}). One can see from Eq.~(\ref{eq:jump_stressfxyz}) that there are $d$ essential interface jump conditions in the present neural network method. On the contrary, the implementation of  IIM in \cite{LIL07} requires 6 interface jump conditions for solving the problem with the dimension $d=2$ as follows.
	\beqs
	-\dbblk{p} \:+\: 2\dbblkL{\mu\frac{\pd\bu}{\pd\bn}}\cdot\bn + \bF\cdot\bn &=& 0, \label{eq:jump_fnt_p}\\
	-\dbblkL{\frac{\pd{p}}{\pd\bn}} \:+\: \dbblk{\bg\cdot\bn} + \frac{\pd}{\pd \bftau}(\bF\cdot\bftau) + 2 \dbblkL{\mu\frac{\pd^2}{\pd \bftau^2}(\bu\cdot\bn)} &=& 0, \label{eq:jump_fnt_dp}\\
	\dbblk{\bu} = {\bf0}, \qquad \dbblk{\mu\div\bu} &=& 0, \label{eq:jump_fnt_u_div}\\
	\dbblkL{\mu\frac{\pd\bu}{\pd\bn}}\cdot\bftau \:+\: \dbblkL{\mu\frac{\pd\bu}{\pd\bftau}}\cdot\bn \:+\: \bF\cdot\bftau &=& 0, \label{eq:jump_fnt_tau1}
	\eeqs
Here, $\bftau$ is the tangential vector along the interface $\Gamma$. For the three-dimensional case, one needs to use even 9 interface jump conditions of velocity and pressure in order to implement the augmented IIM in \cite{WT22}. Thus, it is more advantageous to use the present neural network method to solve the Stokes interface problems than the augmented IIMs regarding the interface jump conditions needed, not to mention the latter ones need to solve the resultant linear systems that lack good symmetric properties.

\section{A discontinuity and cusp capturing PINN}
\label{sec:CuspDCNN}

In this section, we present a discontinuity and cusp capturing physics-informed neural network to solve Stokes equations (\ref{eq:Stokes_jumpform})-(\ref{eq:div1}) with the traction balance equation (\ref{eq:jump_stressfxyz}), subject to the boundary condition (\ref{eq:BCs_D}). To proceed, we assume there exists a smooth level set function $\phi(\bx)$ such that the interior and exterior domains are defined as $\Omega^- = \{\bx \in\Omega\, |\,\phi(\bx)<0\}$ and $\Omega^+ = \{\bx \in\Omega\, |\,\phi(\bx)>0\}$, respectively, with the interface position given by the zero level set, $\Gamma = \{\bx \in\Omega\,|\, \phi(\bx)=0\}$.

Since the pressure has jump discontinuity across the interface $\Gamma$, we follow the idea proposed in \cite{HLL22} and represent the $d$-dimensional discontinuous function $p(\bx)$ by a ($d+1$)-dimensional continuous neural network function $P(\bx,y)$, where $\bx \in \Omega$ and $y \in \mathbb{R}$. We then set the augmented variable $y$ to be the indicator function $I(\bx)$ as
	\beq
	I(\bx) = \left\{\bary{rl} -1, & \bx\in\Omega^-,  \\ 1, & \bx\in\Omega^+. \eary \right.
	\eeq
So the pressure $p(\bx)$ in $\Omega^\pm$ can be represented by the neural network function $P(\bx, y=\pm 1)$. We point out that the neural network $P$ is a smooth extension in $\mathbb{R}^{d+1}$ if a smooth activation function is used. The pressure $p(\bx)$ in $\Omega^-$ and $\Omega^+$ can be regarded as the network function $P$ restricted on two disjoint $d$-dimensional hyperplanes, $S^-=\{(\bx,-1)\:|\:\bx\in\Omega^-\}\subset\mathbb{R}^{d+1}$ and $S^+=\{(\bx,1)\:|\:\bx\in\Omega^+\}\subset\mathbb{R}^{d+1}$, respectively. Consequently, one can rewrite the pressure jump $\dbblk{p}$ as the difference of the continuous network function $P$ evaluated at the interface position on the two disjoint hyperplanes, $S^\pm$. That is, for any $\bx_\Gamma\in\Gamma$,
	\beq
	\dbblk{p}(\bx_\Gamma) = P(\bx_\Gamma,1) - P(\bx_\Gamma,-1)= \dbblk{P}(\bx_\Gamma, 0). \label{eq:p_jump_netw}
	\eeq		
The gradient of $p$ in $\Omega^\pm$ can be written by
	\beq
	\grad p(\bx) = \grad_\bx P(\bx,I(\bx)) + \pd_y P(\bx,I(\bx)) \, \grad I(\bx) = \grad_\bx P(\bx,I(\bx)),
	\eeq
	where $\nabla_\bx P\in\mathbb{R}^d$ represents a vector consisting of the partial derivatives of $P$ with respect to the components in $\bx$, and $\pd_y P$ is the shorthand notation of partial derivative of $P$ with respect to the augmented variable $y$.

Similarly, following the idea in \cite{TLHL22},  we represent the $d$-dimensional non-smooth velocity function $\bu(\bx)$ by a ($d+1$)-dimensional smooth neural network function $\bU(\bx,z)=(U_1(\bx, z), \ldots U_d(\bx, z))$.  We then set the augmented variable $z$ to be $z=\phi_a(\bx)=\abs{\phi(\bx)}$ so the velocity $\bu(\bx)$ in $\Omega$ can be represented by a neural network $\bU(\bx, z=\phi_a(\bx))$. Here, we call $\phi_a(\bx)$ as a cusp-enforced level set function since it exhibits a gradient jump defined as $\dbblk{\grad \phi_a}(\bx_\Gamma) = 2 \grad \phi(\bx_\Gamma)$ for $\bx_\Gamma \in \Gamma$. Based on the above property, one can immediately see that the gradient of velocity $\grad u_k = \grad_\bx U_k + \pd_z U_k \grad \phi_a, k=1,\ldots,d$ in $\Omega^\pm$ has discontinuity jump across the interface which inherits the cusp-like solution behavior of the velocity. More precisely, the jump of $\mu\grad u_k$ across the interface can be derived directly as
	\beqs
	\dbblk{\mu\grad u_k}(\bx_\Gamma) &=& \dbblk{\mu\grad_\bx U_k}(\bx_\Gamma) + \dbblk{\mu\pd_z U_k \grad \phi_a}(\bx_\Gamma) \nonumber\\
		&=& \dbblk{\mu}\grad_\bx U_k(\bx_\Gamma) + (\mu^++\mu^-)\pd_z U_k(\bx_\Gamma)\grad \phi(\bx_\Gamma), \label{eq:term1}
	\eeqs
for $k=1,\ldots,d$. Similarly, the jump of $\mu\frac{\pd\bu}{\pd{x}_k}$ can be calculated by
\beqs
	\dbblkL{\mu\frac{\pd\bu}{\pd{x}_k}}(\bx_\Gamma) &=& \dbblkL{\mu\frac{\pd\bU}{\pd{x}_k}}(\bx_\Gamma) + \dbblkL{\mu
\frac{\pd\bU}{\pd{z}}\frac{\pd\phi_a}{\pd{x}_k} }(\bx_\Gamma) \nonumber\\
		&=& \dbblk{\mu}\frac{\pd\bU}{\pd{x}_k}(\bx_\Gamma) + (\mu^++\mu^-) \frac{\pd\bU}{\pd{z}}\frac{\pd\phi}{\pd{x}_k}(\bx_\Gamma). \label{eq:term2}
	\eeqs
By multiplying the normal vector $\bn=\grad \phi/ \|\grad \phi \|$ to the above equations (\ref{eq:term1}) and (\ref{eq:term2}), respectively; we obtain the following derivative jumps of $\bu$ as follows.
	\beqs
 	\dbblkL{\mu\frac{\pd u_k}{\pd \bn}}(\bx_\Gamma) &=& \dbblk{\mu}\left(\grad_\bx U_k \cdot\bn\right)(\bx_\Gamma) + (\mu^++\mu^-)\frac{\pd U_k}{\pd{z}}\|\grad\phi\|(\bx_\Gamma), \label{eq:stress_d}\\
	\dbblkL{\mu\frac{\pd\bu}{\pd{x}_k}}\cdot\bn(\bx_\Gamma) &=& \dbblk{\mu}\left(\frac{\pd\bU}{\pd{x}_k}\cdot\bn\right)(\bx_\Gamma) + (\mu^++\mu^-) \left(\frac{\pd\bU}{\pd{z}}\cdot\bn\right)\frac{\pd\phi}{\pd{x}_k}(\bx_\Gamma), \label{eq:stress_t}
	\eeqs
for $k=1,\ldots,d$. Therefore, the interfacial velocity derivative jump terms in Eq.~(\ref{eq:jump_stressfxyz}) can be expressed by the smooth network function $\bU(\bx, z)$ with the assistance of the level set function $\phi$.

After careful calculation, the Laplacian of $u_k$ in $\Omega^\pm$ can also be explicitly expressed in terms of $U_k$ and $\phi_a$ as
	\beq
	\label{eq:lapU}
	\mu\Delta u_k = \mu\left(\Delta_\bx U_k + 2\grad\phi_a\cdot\grad_\bx\left( \pd_{z}U_k\right)+\|\grad\phi_a\|^2\pd_{zz}U_k+\pd_z U_k\Delta\phi_a\right),
	\eeq
	where $\Delta_\bx$ represents the Laplace operator applied only to the variable $\bx$.
	Since the neural network function $\bU$ is smooth, the continuity of velocity vector holds automatically, and calculating the derivatives of the network $\bU$ with respect to its input variables $\bx$ and $z$ via automatic differentiation~\cite{GW08} is straightforward.

Next, we express the Stokes equations~(\ref{eq:Stokes_jumpform})-(\ref{eq:div1}) in terms of $P$ and $\bU$ as follows. To simplify, we introduce the notation $\mathcal{L} U_k$ to represent the right-hand side of  $\mu\Delta u_k$ in Eq.~(\ref{eq:lapU}) so that Eqs.~(\ref{eq:Stokes_jumpform}) and~(\ref{eq:div1}) are written as
	\beqs
	-\frac{\pd{P}}{\pd{x_k}}(\bx,I(\bx))+\mathcal{L} U_k(\bx,\phi_a(\bx))+ g_k(\bx) &=& 0, \quad \bx \in \Omega^\pm,\quad k=1,\ldots,d,\label{eq:Stokes_netPU}\\
	\grad_\bx\cdot\bU(\bx,\phi_a(\bx)) + \frac{\pd\bU}{\pd{z}}(\bx,\phi_a(\bx))\cdot\grad\phi_a(\bx) &=& 0, \quad \bx \in \Omega^\pm. \label{eq:divU}
	\eeqs
We also denote the sum of the right-hand sides of Eqs.~(\ref{eq:stress_d})-(\ref{eq:stress_t}) by $\mathcal{D}_{k}\bU$, so the traction balance equation (\ref{eq:jump_stressfxyz}) has the form as
	\beq
	-\dbblk{P}(\bx_\Gamma,0)n_k(\bx_\Gamma) + \mathcal{D}_{k}\bU(\bx_\Gamma,0) + F_k(\bx_\Gamma) = 0, \quad \bx_\Gamma \in \Gamma,\quad k=1,\ldots,d. \label{eq:jump_stress_PU}
	\eeq
Regarding the boundary condition Eq.(\ref{eq:BCs_D}), we use $\bx_B$ to represent the points located on $\pd\Omega$ so it becomes
	\beq
	U_k(\bx_B,\phi_a(\bx_B)) - u_{b,k}(\bx_B) = 0, \quad \bx_B\in\pd\Omega, \quad k=1,\ldots,d. \label{eq:BCs_D_U}
	\eeq
	
\begin{figure}[h]
	\begin{center}
	\includegraphics[scale=0.4]{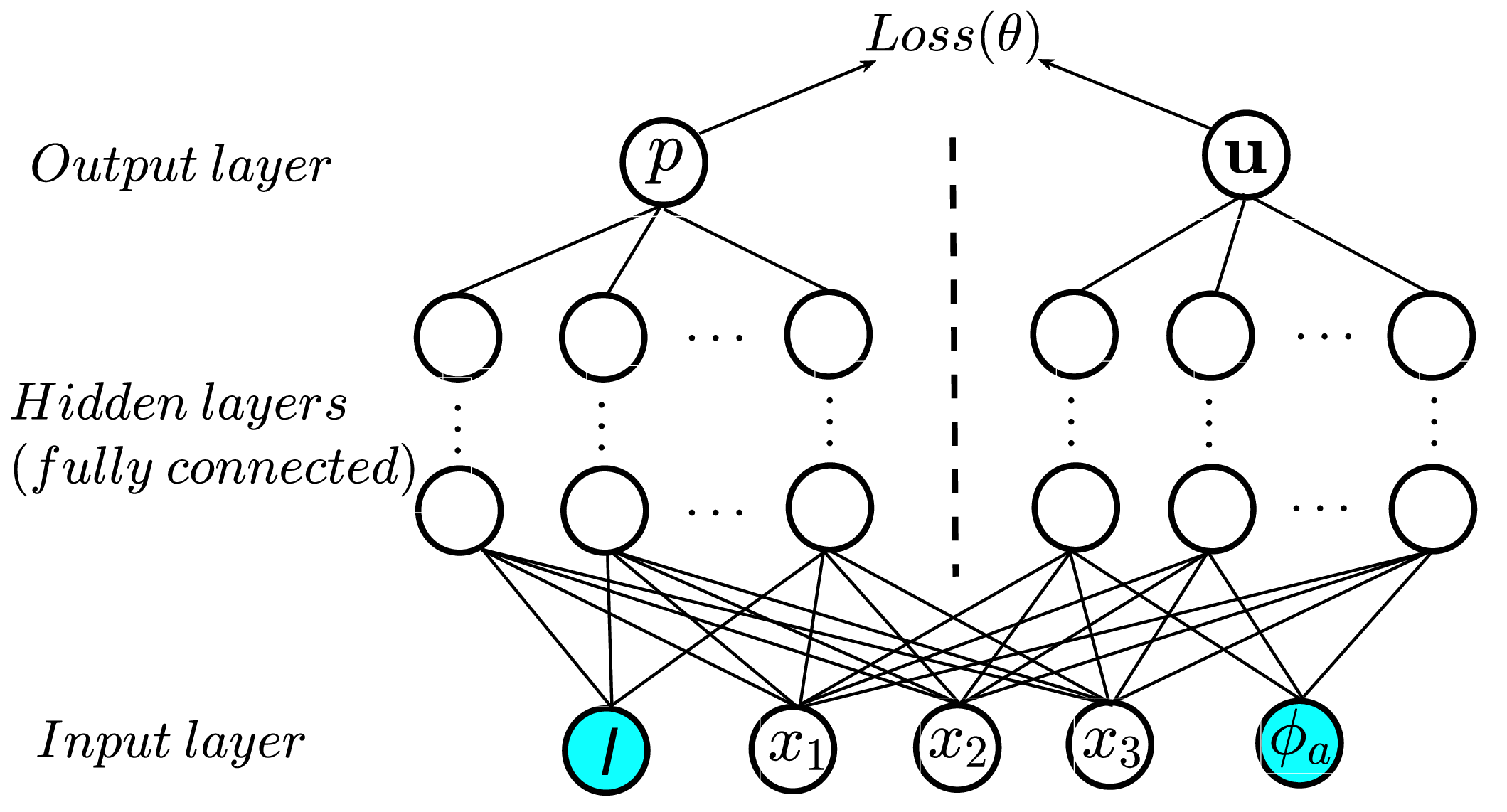}
	\caption{Diagram of the $L$-hidden-layer network structure.}\label{fig:DC-CuspNN}
	\end{center}
	\end{figure}

Next, we aim to show the network architecture for the pressure $P(\bx, I(\bx))$ and the velocity $\bU(\bx, \phi_a(\bx))$.  Figure~\ref{fig:DC-CuspNN} illustrates the present neural network structure, where $(\bx,I(\bx))\in\mathbb{R}^{d+1}$ represents the $d+1$ feature input for the sub-network $P$, while $(\bx,\phi_a(\bx))\in\mathbb{R}^{d+1}$ represents the input for the sub-network $\bU$. Although the two sub-networks have separate structures, their output layers are fed into the same loss function, $\Loss(\bft)$, which will be given in detail later. For simplicity, we assume that both sub-networks have the same depth number, $L$, and we use the subscripts $\star=$'$p$', '$u$' to distinguish which sub-network is present. In each sub-network, we label the input layer as layer $0$ and denote the feature input as $\bv_\star^{[0]}=(\bx,I(\bx))^T$ or $(\bx,\phi_a(\bx))^T$, a column vector in $\mathbb{R}^{d+1}$. The output at the $\ell$-th hidden layer, consisting of $N_\star^{[\ell]}$ neurons, denoted as $\bv_\star^{[\ell]}\in\mathbb{R}^{N_\star^{[\ell]}}$, represents an affine mapping of the output of layer $\ell-1$ (i.e., $\bv_\star^{[\ell-1]}$), followed by an element-wise activation function $\sigma$ as
	\beq
	\bv_\star^{[\ell]} = \sigma\left(W_\star^{[\ell]}\bv_\star^{[\ell-1]}+\bb_\star^{[\ell]}\right), \quad \ell=1,\cdots,L. \label{eq:DNN}
	\eeq
Here, $W_\star^{[\ell]}\in\mathbb{R}^{N_\star^{[\ell]}\times N_\star^{[\ell-1]}}$ represents the weight matrix connecting layer $\ell-1$ to layer $\ell$, and $\bb_\star^{[\ell]}\in\mathbb{R}^{N_\star^{[\ell]}}$ is the bias vector at layer $\ell$. Finally, the output of this $L$-hidden-layer network is denoted by the following
	\beqs
	P_{\mathcal{N}}(\bx, I(\bx);\bft_p) &=& W_p^{[L+1]}\mathbf{v}_p^{[L]}, \\
	\bU_{\mathcal{N}}(\bx, \phi_a(\bx);\bft_u) &=& W_u^{[L+1]}\mathbf{v}_u^{[L]},
	\eeqs
where $W_p^{[L+1]}\in\mathbb{R}^{1\times N_p^{[L]}}$ and $W_u^{[L+1]}\in\mathbb{R}^{d\times N_u^{[L]}}$. The notation $\bft_\star$ denotes the vector collecting all trainable parameters (including all the weights and biases) in each corresponding sub-network. We use $\bft=(\bft_p,\bft_u)$ to represent the collection of all trainable parameters within the whole network, and the dimension of $\bft$ is counted as
	\beq
	N_\bft = \left( N_p^{[L]}+\sum^{L}_{\ell=1}(N_p^{[\ell-1]}+1)N_p^{[\ell]} \right) + \left( N_u^{[L]}d+\sum^{L}_{\ell=1}(N_u^{[\ell-1]}+1)N_u^{[\ell]} \right).
	\eeq
	
In the training process, we randomly choose $M_I$ points, $\left\{\bx^i\right\}_{i=1}^{M_I}$, within the region of $\Omega^-\cup\Omega^+$, and $M_B$ points, $\left\{\bx^i_B\right\}_{i=1}^{M_B}$, on the boundary of the domain $\partial\Omega$, respectively. Additionally, we select $M_\Gamma$ points, $\left\{\bx^i_\Gamma\right\}_{i=1}^{M_\Gamma}$, on the interface $\Gamma$. Altogether, we have a total number of $M=M_I+M_B+M_\Gamma$ training points. Within the PINN learning framework, the loss function consists of the mean squared errors of the residual  of differential equations~(\ref{eq:Stokes_netPU})-(\ref{eq:divU}), the traction balance equation (\ref{eq:jump_stress_PU}), and the boundary condition (\ref{eq:BCs_D_U}).	That is,
\beqs
\label{eq:loss}
\Loss(\bft) &=& \frac{1}{M_{I}}\sum_{i=1}^{M_{I}} \left(\sum^d_{k=1}\abs{L_{I,k}(\bx^i,I(\bx^i),\phi_a(\bx^i);\bft)}^{2} + \abs{L_D(\bx^i,\phi_a(\bx^i);\bft_u)}^{2} \right)   \\
	&+& \frac{c_\Gamma}{M_\Gamma}\sum_{i=1}^{M_\Gamma} \sum^d_{k=1} \abs{ L_{\Gamma,k}(\bx^i_\Gamma,0;\bft)}^{2} +
	\frac{c_B}{M_{B}}\sum_{i=1}^{M_{B}} \sum^d_{k=1}\abs{L_{B,k}(\bx^i_B,\phi_a(\bx^i_B);\bft_u)}^{2}. \nonumber
\eeqs
where the residual error $L_I$, incompressibility constraint error $L_D$, traction balance condition error $L_\Gamma$, and boundary condition error $L_B$, are shown respectively as follows:
\begin{eqnarray}
L_{I,k}(\bx,I(\bx),\phi_a(\bx);\bft) &=& -\frac{\pd{P}_\mathcal{N}}{\pd{x_k}}(\bx,I(\bx);\bft_p)+\mathcal{L} U_{\mathcal{N},k}(\bx,\phi_a(\bx);\bft_u)-g_k(\bx), \label{eq:lossi}\\
L_D(\bx,\phi_a(\bx);\bft_u) &=& \grad_\bx\cdot\bU_\mathcal{N}(\bx,\phi_a(\bx);\bft_u) + \frac{\pd\bU_\mathcal{N}}{\pd{z}}(\bx,\phi_a(\bx);\bft_u)\cdot\grad\phi_a(\bx), \label{eq:lossd}\\
L_{\Gamma,k}(\bx_\Gamma,0;\bft) &=& -\dbblk{P_\mathcal{N}}(\bx_\Gamma,0;\bft_p)n_k(\bx_\Gamma) + \mathcal{D}_{k}\bU_\mathcal{N}(\bx_\Gamma,0;\bft_u)+ F_k(\bx_\Gamma),\qquad \label{eq:lossgk}\\
L_{B,k}(\bx_B,\phi_a(\bx_B);\bft_u) &=& U_{\mathcal{N},k}(\bx_B,\phi_a(\bx_B);\bft_u)-u_{b,k}(\bx_B). \label{eq:lossb}
\end{eqnarray}
The constants $c_B$ and $c_\Gamma$ appeared in the loss function~(\ref{eq:loss}) are chosen to balance the contribution of the terms related to the traction balance equation ~(\ref{eq:jump_stress_PU}), and the boundary condition~(\ref{eq:BCs_D_U}), respectively.

\section{Numerical results}
\label{sec:numerics}
In this section, we perform a series of accuracy tests for the proposed discontinuity and cusp capturing neural network method on solving two- and three-dimensional Stokes interface problems described by Eqs.~(\ref{eq:Stokes_jumpform})-(\ref{eq:BCs_D}). We set the penalty constants in the loss function, $c_\Gamma=c_B=1$. This approach enables the use of smooth neural network functions defined in $\mathbb{R}^{d+1}$, denoted by $P(\bx,y)$ and $\bU(\bx,z)$, to learn discontinuous and non-smooth solutions in $\mathbb{R}^{d}$, $p(\bx)$ and $\bu(\bx)$, respectively. We denote the network-predicted solutions as $P_\mathcal{N}$ and $\bU_\mathcal{N}=(U_{\mathcal{N},1},\ldots,U_{\mathcal{N},d})$. Throughout this section, we represent the interface $\Gamma$ using the zero level set of a suitable level set function $\phi(\bx)$. We use the indicator function $I(\bx)$ and the cusp-enforced level set function $\phi_a(\bx)=|\phi(\bx)|$ as the augmented input for the network of $P_\mathcal{N}$ and $\bU_\mathcal{N}$, respectively.

To ensure the $C^2$-regularity of $\bu(\bx)$ in each subdomain, we employ the sigmoid function $\sigma(x)=\frac{1}{1+e^{-x}}$ as the activation function. For the numerical examples presented here (except the last example whose solution is not available, so deeper neural network approximations are used), we adopt a completely shallow network structure ($L=1$) with $N_p$ and $N_u$ neurons for the pressure and velocity, respectively, as shown in the network architecture Figure~\ref{fig:DC-CuspNN}. So the number of total hyper-parameters to be trained is $N_\bft = (d+3)N_p + 2(d+1)N_u$. The training and test data points are generated using the Latin hypercube sampling algorithm~\cite{MBC79}. This algorithm effectively avoids clustering data points at specific locations, resulting in a nearly random sampling. To assess the accuracy of the network solution, we select $M_{test}$ points (distinct from the training points) within $\Omega$ and calculate the $L^\infty$ errors as follows:
	\beqs
	E_p^\infty = \|p-P_\mathcal{N}\|_\infty, \quad E_\bu^\infty = \frac{1}{d}\sum^d_{k=1}\|u_k-U_{\mathcal{N},k}\|_\infty.\nonumber
	\eeqs
	Here, the $L^\infty$ error of a function $f$ is defined as $\|f\|_\infty = \max_{1\leq i \leq M_{test}}|f(\bx^i)|$. We set the number of test points for each experiment to be $M_{test}=100M$, where $M$ represents the total number of training points used. Since the predicted results will vary slightly due to the randomness of training and test data points and the initialization of trainable parameters, we report the average value of errors and losses over 5 trial runs.

During the training process, we employ the Levenberg-Marquardt (LM) algorithm~\cite{More78} as the optimizer to train the network and update the damping parameter $\nu$ using the strategies outlined in~\cite{TS12}. The training process is terminated either when the loss value $\Loss(\bft)$ falls below a threshold $\epsilon_\theta=10^{-14}$  or the maximum iteration step $Epoch_{max}=3000$ is reached, unless otherwise stated. All experiments are performed on a desktop computer equipped with a single NVIDIA GeForce RTX3060 GPU. We implement the present network architecture using PyTorch (v1.13)~\cite{PGMetel19}, and all trainable parameters (i.e., weights and biases) are initialized using PyTorch's default settings. We have made the source code used in the paper available on GitHub at "https://github.com/yuhautseng/".	

Let us remark on the choice of test examples in the following. As mentioned earlier, we aim to emphasize the ease of code implementation for the problem compared with the traditional IIM; thus, we have taken the test examples in 2D and 3D directly from the literature \cite{LIL07} and \cite{WT22}, respectively. One might wonder why the interfaces $\Gamma$ chosen there are in circular (2D) or spherical (3D) form. The reason is quite apparent since it is easy to construct the analytic solution under such a circumstance. (The solution must satisfy the jump conditions (\ref{eq:jump_stressfxyz}) for $k=1, \ldots d$.) Nevertheless, this kind of simple interface shape cannot be regarded as a limitation of the present neural network methodology since it is completely mesh-free. Indeed, it can be applied to more complex interface cases when the analytic solutions are not available (see Example 5).

Notice that, it is not completely fair to make a running time comparison between the present neural network method  with IIM since the former method is mesh-free while the latter one is grid-based. The present PINN method learns the solution directly while the IIM solves the solution on uniform grids. Besides, the results in \cite{LIL07} and \cite{WT22} were run in different computing hardware and it is unlikely for us to reproduce the results. To provide some idea of the time performance of the present PINN method, we rerun our code in a slightly faster desktop machine equipped with a single NVIDIA GeForce RTX4090 GPU. The total training time up to $3000$ epochs for the two-dimensional case with $(N_p, M_0)=(20,30)$ (as shown in second row of Table 1) cost approximately $59$ seconds while the three-dimensional case with $(N_p, N_u)=(40,100)$ (as shown in third row of Table 3) cost approximately $253$ seconds, respectively. \\
		
\noindent{\bf Example 1}: We start by solving a two-dimensional Stokes interface problem with a piecewise-constant viscosity and a singular force in a domain $\Omega=[-2,2]^2$. This example is taken directly from Example 4.2 in \cite{LIL07} so that we can compare the results with the augmented immersed interface method developed there. The interface $\Gamma$ is simply a unit circle defined as the zero level set of the function $\phi(\bx) = x_1^2+x_2^2-1$. We specify the viscosity $\mu(\bx)$, and the exact solutions $p(\bx)$ and $\bu(\bx)=(u_1(\bx),u_2(\bx))$ as
\begin{align*}
\mu(\bx) &= \begin{cases} 1, & \text{in }\Omega^-, \\ 0.5, & \text{in }\Omega^+, \end{cases} &
p(\bx)   &= \begin{cases} 1, & \text{in }\Omega^-, \\ 0,   & \text{in }\Omega^+, \end{cases} \\
u_1(\bx) &= \begin{cases} x_2(x_1^2+x_2^2-1),  & \text{in }\Omega^-, \\ 0, & \text{in }\Omega^+, \end{cases} &
u_2(\bx) &= \begin{cases} -x_1(x_1^2+x_2^2-1), & \text{in }\Omega^-, \\ 0, & \text{in }\Omega^+. \end{cases}
\end{align*}
The external force $\bg(\bx)=(g_1(\bx),g_2(\bx))$ thus can be derived accordingly as
\begin{align*}
g_1(\bx) &= \begin{cases} -8x_2, & \text{in }\Omega^-, \\ 0, & \text{in }\Omega^+, \end{cases} &
g_2(\bx) &= \begin{cases} 8x_1, & \text{in }\Omega^-, \\ 0, & \text{in }\Omega^+. \end{cases}
\end{align*}
and the interfacial force $\bF(\bx_\Gamma)=(F_1(\bx_\Gamma),F_2(\bx_\Gamma))=(-x_1-2x_2, -x_2+2x_1)$ can also be obtained from Eq.~(\ref{eq:jump_stressfxyz}).

As mentioned earlier, we use a completely shallow network structure ($L=1$) with $N_p$ and $N_u$ neurons for the pressure $P_\mathcal{N}$ and velocity $\bU_\mathcal{N}$, respectively. Here, we choose $N_p=10,\,20,\, 30$, and $N_u=2N_p$ so that the total training parameters is $N_\bft=17N_p$. We also vary the training points accordingly by introducing a grid parameter $M_0=20,\,30,\,40$ (can be regarded as the grid number used in traditional grid-based methods) so that the training dataset consists of $M_I = M_0^2$ points in $\Omega^-\cup\Omega^+$, $M_\Gamma = 3M_0$ points on the interface $\Gamma$, and $M_B = 4M_0$ points on the boundary $\partial\Omega$. The total number of training points is $M=M_0(M_0+7)$.

Table~\ref{tab:ex1_Fxy_vs_IIM} shows the $L^\infty$ errors of $P_\mathcal{N}$ and $\bU_\mathcal{N}$, as well as the corresponding training losses. We also show the results obtained by IIM in \cite{LIL07} using various grid resolutions in the right panel for comparison. It can be observed that the proposed network method achieves prediction accuracy on the order of $O(10^{-5})$ to $O(10^{-8})$ as $(N_p, M_0)$ varies from $(10, 20)$ to $(30,40)$, with the loss dropping from $O(10^{-10})$ to $O(10^{-15})$ accordingly. Meanwhile, Table~\ref{tab:ex1_Fxy_vs_IIM} can be considered as an informal convergence analysis associated with the number of learnable parameters and training points used in the network. As seen, using a resolution of $512^2$ uniform mesh, the IIM yields the accuracy $E_p^\infty$ and $E_\bu^\infty$ of the magnitude $O(10^{-5})$ that is at least two orders of magnitude greater than the one obtained by the present network using relatively fewer training points $M=1880$. \\
\begin{table}[hbt!]
\begin{center}
\begin{tabular}{cccc|ccc}
\hline
& & Present & & & IIM~\cite{LIL07} & \\
\hline
$(N_p, M_0)$ & $E_p^\infty$ & $E_\bu^\infty$ & $\Loss(\bft)$ & Grid No.& $E_p^\infty$ & $E_\bu^\infty$\\
\hline
$(10,20)$ & $8.94\times10^{-5}$ & $1.16\times10^{-5}$ & $6.55\times10^{-10}$ & $128^2$ & $8.10\times10^{-4}$ & $2.27\times10^{-4}$ \\
$(20,30)$ & $1.50\times10^{-6}$ & $2.73\times10^{-7}$ & $2.54\times10^{-14}$ & $256^2$ & $2.54\times10^{-4}$ & $4.77\times10^{-5}$ \\
$(30, 40)$ & $4.05\times10^{-7}$ & $6.87\times10^{-8}$ & $9.04\times10^{-15}$ & $512^2$ & $1.41\times10^{-5}$ & $1.41\times10^{-5}$ \\
\hline
\end{tabular}
\end{center}
\caption{The $L^\infty$ errors of $P_\mathcal{N}$ and $\bU_\mathcal{N}$ for Example 1. Left panel: the present network result. Here, $N_\bft=17N_p$, $M=M_0(M_0+7)$, $\epsilon_\theta=10^{-14}$, and $Epoch_{max}=3000$. Right panel: IIM results in \cite{LIL07} with grid resolutions $128^2$, $256^2$, and $512^2$.} \label{tab:ex1_Fxy_vs_IIM}
\end{table}

\noindent{\bf Example 2}: In the second example (also referred to as Example 4.3 in~\cite{LIL07}), we keep the domain $\Omega$ and interface $\Gamma$ same as shown in Example 1. However, we choose the pressure jump no longer as a constant and it is defined by
\begin{align}
p(\bx) &= \begin{cases} \left(\frac{3}{8}-\frac{3}{4}x_1^2\right)x_1x_2, & \text{in } \Omega^-, \\ 0, & \text{in } \Omega^+. \end{cases} \nonumber
\end{align}
The exact solution of velocity $\bu(\bx)=(u_1(\bx),u_2(\bx))$ is chosen as
\begin{align}
u_1(\bx) &= \begin{cases} \frac{1}{4}x_2, & \text{in } \Omega^-, \\ \frac{1}{4}x_2(x_1^2+x_2^2), & \text{in } \Omega^+, \end{cases} &
u_2(\bx) &= \begin{cases} -\frac{1}{4}x_1(1-x_1^2), & \text{in } \Omega^-, \\ -\frac{1}{4}x_1x_2^2, & \text{in } \Omega^+, \end{cases} \nonumber
\end{align}
and the external force $\bg(\bx)=(g_1(\bx),g_2(\bx))$ is derived as
\begin{align}
g_1(\bx) &= \begin{cases} \left(\frac{3}{8}-\frac{9}{4}x_1^2\right)x_2, & \text{in } \Omega^-,\\ -2\mu^+x_2, & \text{in } \Omega^+, \end{cases} &
g_2(\bx) &= \begin{cases} \left(\frac{3}{8}-\frac{3}{2}\mu^--\frac{3}{4}x_1^2\right)x_1, & \text{in } \Omega^-,\\ \frac{1}{2}\mu^+x_1, & \text{in } \Omega^+, \end{cases} \nonumber
\end{align}
The interfacial force, $\bF(\bx_\Gamma)=(F_1(\bx_\Gamma),F_2(\bx_\Gamma))$, can be derived from the jump condition~(\ref{eq:jump_stressfxyz}) as
\begin{align}
F_1(\bx_\Gamma) &= \left(\frac{3}{4}x_1^2-\frac{3}{8}\right)x_1^2x_2 - \frac{3}{2}\dbblk{\mu}x_1^4x_2 + \frac{1}{2}\mu^+x_2 + \frac{3}{4}\dbblk{\mu}x_1^2(1-2x_1^2)x_2, \nonumber \\
F_2(\bx_\Gamma) &= \left(\frac{3}{4}x_1^2-\frac{3}{8}\right)x_1x_2^2 - \frac{3}{2}\dbblk{\mu}x_1^3x_2^2 - \frac{1}{2}\mu^+x_1 - \frac{3}{4}\dbblk{\mu}x_1^3(1-2x_1^2). \nonumber
\end{align}

For the network structure, we choose the exactly same setup as in Example 1. That is, we use the same values of $N_p$ and $M_0$ but vary the viscosity contrast as $(\mu^-,\mu^+)=(1,0.1), (10^{-3}, 1), (1, 10^{-3})$.  Table~\ref{tab:ex2_Fxy_vs_IIM} shows the $L^\infty$ norm errors of $P_\mathcal{N}$ and $\bU_\mathcal{N}$, and the corresponding training losses in the left panel. The corresponding results from IIM~\cite{LIL07} are shown in the right panel for comparison. In the first case $(\mu^-,\mu^+)=(1,0.1)$, by varying $(N_p, M_0)$ from $(10, 20)$ to $(30, 40)$, we again observe that the present network predictions $P_\mathcal{N}$ and $\bU_\mathcal{N}$ can achieve the accuracy with $L^\infty$ errors ranging from the order $O(10^{-5})$ to $O(10^{-7})$ and the loss drops from $O(10^{-10})$ to $O(10^{-14})$ accordingly. As we increase the neurons $N_p$ and the training points $M_0$, we again see an informal evidence for the numerical convergence of the present method.
The results in~\cite{LIL07} have $E_p^\infty=1.54\times10^{-4}$ and $E_\bu^\infty=6.49\times10^{-5}$ under the resolution using a $512\times512$ uniform mesh. Meanwhile, we also show the capability of the present method for interface problems with high viscosity contrast as $(\mu^-,\mu^+)=(10^{-3},1)$ and $(1,10^{-3})$. The present method achieves accuracy with $L^\infty$ errors ranging from $E^\infty_p=O(10^{-4})$ to $E^\infty_\bu=O(10^{-7})$ which outperforms the results obtained from IIM significantly. \\
	\begin{table}[hbt!]
	\begin{center}
   	\begin{tabular}{c|cccc|ccc}
   	\hline
   		&	& & Present & & & IIM~\cite{LIL07} & \\
   	\hline
   	$(\mu^-,\mu^+)$ &	$(N_p, M_0)$   & $E_p^\infty$ & $E_\bu^\infty$ & $\Loss(\bft)$ & $N_{IIM}$ & $E_p^\infty$ & $E_\bu^\infty$ \\
   	\hline
   	      & $(10, 20)$ & $4.43\times10^{-5}$ & $7.43\times10^{-6}$ & $2.83\times10^{-10}$ & $128^2$ & $2.30\times10^{-3}$ & $1.21\times10^{-3}$ \\
$(1,0.1)$ & $(20, 30)$ & $3.14\times10^{-6}$ & $5.49\times10^{-7}$ & $1.34\times10^{-12}$ & $256^2$ & $5.47\times10^{-4}$ & $2.69\times10^{-4}$ \\
	      & $(30, 40)$ & $1.08\times10^{-6}$ & $1.21\times10^{-7}$ & $5.09\times10^{-14}$ & $512^2$ & $1.54\times10^{-4}$ & $6.49\times10^{-5}$ \\
   \hline
   &  $(10, 20)$ & $5.64\times10^{-4}$ & $8.74\times10^{-5}$ & $3.00\times10^{-9}$ & $128^2$ & $1.04\times10^{-3}$ & $6.23\times10^{-2}$ \\
$(10^{-3},1)$ & $(20, 30)$ & $5.84\times10^{-5}$ & $2.59\times10^{-6}$ & $6.84\times10^{-12}$ & $256^2$ & $3.59\times10^{-4}$ & $1.40\times10^{-2}$ \\
	          & $(30, 40)$ & $2.65\times10^{-6}$ & $2.23\times10^{-7}$ & $4.62\times10^{-14}$ & $512^2$ & $7.09\times10^{-5}$ & $2.82\times10^{-3}$ \\
   \hline
              & $(10, 20)$& $6.62\times10^{-4}$ & $1.15\times10^{-4}$ & $2.42\times10^{-8}$ & $128^2$ & $6.53\times10^{-3}$ & $3.15\times10^{-1}$ \\
$(1,10^{-3})$ & $(20, 30)$ & $5.78\times10^{-5}$ & $6.28\times10^{-6}$ & $1.98\times10^{-10}$ & $256^2$ & $1.18\times10^{-3}$ & $4.64\times10^{-2}$ \\
	          & $(30, 40)$ & $1.10\times10^{-6}$ & $1.44\times10^{-7}$ & $6.72\times10^{-14}$ & $512^2$ & $3.02\times10^{-4}$ & $1.17\times10^{-3}$ \\
   \hline
   \end{tabular}
   \end{center}
   \caption{The $L^\infty$  errors of $P_\mathcal{N}$ and $\bU_\mathcal{N}$ for Example 2. Left panel:  the present network result. Here, $N_\bft=17N_p$, $M=M_0(M_0+7)$, $\epsilon_\theta=10^{-14}$, and $Epoch_{max}=3000$.  Right panel: IIM  results in \cite{LIL07} with grid resolutions $128^2$, $256^2$ and $512^2$.} \label{tab:ex2_Fxy_vs_IIM}	
\end{table}

\noindent{\bf Example 3}: We now consider solving a three-dimensional Stokes interface problem with a piecewise-constant viscosity ($\mu^-=1$ and $\mu^+=0.1$) and singular force in a domain $\Omega=[-2,2]^3$. The interface $\Gamma$ is a simple unit sphere so it can be described by the zero level set of the function $\phi(\bx) = \|\bx\|^2_2-1$.  This example is taken directly from Example 6.2 in \cite{WT22} so that we can compare the results with the augmented immersed interface method developed there. We list the exact solutions of pressure $p(\bx)$ and velocity  $\bu(\bx)$ as follows:
\beq
p(\bx) = \left\{ \bary{ll} 1, & \Omega^-,\\ 0, & \Omega^+, \eary\right. \qquad
\bu(\bx)=\left(\bary{c}u_1(\bx) \\ u_2(\bx) \\ u_3(\bx) \eary\right) =
	\left(\bary{r} x_2x_3\left(\|\bx\|^2_2-1\right) \\ x_3x_1\left(\|\bx\|^2_2-1\right) \\ -2x_1x_2\left(\|\bx\|^2_2-1\right) \eary\right).
	\nonumber
\eeq
Here, the velocity components $u_1(\bx)$, $u_2(\bx)$, and $u_3(\bx)$ are all smooth functions over the entire $\Omega$, and $p(\bx)$ is a piecewise constant. Accordingly, the external source $\bg(\bx)=(g_1(\bx),g_2(\bx),g_3(\bx))$ can be derived as
\beq
g_1(\bx) = \left\{ \bary{ll} -14\mu^-x_2x_3, & \Omega^-,\\ -14\mu^+x_2x_3, & \Omega^+, \eary\right. \:\:
g_2(\bx) = \left\{ \bary{ll} -14\mu^-x_3x_1, & \Omega^-,\\ -14\mu^+x_3x_1, & \Omega^+, \eary\right.  \:\:
g_3(\bx) = \left\{ \bary{ll} 28\mu^-x_1x_2, & \Omega^-,\\ 28\mu^+x_1x_2, & \Omega^+. \eary\right.  \nonumber
\eeq
Again, the interfacial force can be derived from the jump condition~(\ref{eq:jump_stressfxyz}) as
\beq
\bF(\bx) = \left(\bary{c} F_1(\bx) \\ F_2(\bx) \\ F_3(\bx) \eary\right) = \left(\bary{c}-x_1 - 2\dbblk{\mu}x_2x_3 \\ -x_2 - 2\dbblk{\mu}x_3x_1 \\ -x_3 + 4\dbblk{\mu}x_1x_2\eary\right). \nonumber
\eeq

	\begin{table}[h!]
	\begin{center}
   	\begin{tabular}{cccc|cc}
   	\hline
     & & Present & &  &  IIM~\cite{WT22} \\
   	\hline
   	$(N_p,N_u)$ & $E^\infty_p$ & $E^\infty_\bu$ & $\Loss(\bft)$ & Grid no. &  $E^\infty_p$\qquad\qquad$E^\infty_\bu$ \\
   	\hline
   	$(25,60)$   & $5.04\times10^{-3}$ & $2.61\times10^{-3}$ & $8.73\times10^{-7}$ & $64^3$  &  $3.06\times10^{-3}$\quad$2.54\times10^{-3}$ \\
   	$(30,75)$   & $5.02\times10^{-4}$ & $1.91\times10^{-4}$ & $6.80\times10^{-9}$ & $128^3$ &  $6.44\times10^{-4}$\quad$5.40\times10^{-4}$ \\
   	$(40,100)$ & $5.06\times10^{-5}$ & $3.92\times10^{-5}$ & $7.01\times10^{-11}$ & $256^3$ &  $1.68\times10^{-4}$\quad$1.22\times10^{-4}$ \\
 	\hline
    \end{tabular}
	\caption{The $L^\infty$ errors of $P_\mathcal{N}$ and $\bU_\mathcal{N}$ for Example 3. Left panel:  the present network results. Here, $\mu^-=1$ and $\mu^+=0.1$, $M=6216$ training data points. Right panel: IIM  results in  \cite{WT22} with grid resolutions $64^3$, $128^3$, and $256^3$. }
	\label{tab:ex3_error_compare}
    \end{center}	
	\end{table}

Unlike the two-dimensional cases in previous examples, we now fix the number of  training points but vary the number of neurons $N_p$ and $N_u$ to ease our computations. For the sampling of training data points, we generate $M_I$ data points in the region $\Omega^-\cup\Omega^+$, $M_B$ on the domain boundary ($M_B/6$ on each face of $\pd\Omega$), and $M_\Gamma$ data points on the surface $\Gamma$ generated using DistMesh~\cite{PS04}. The total number of training points used in each of the following numerical runs is $M=6216$ ($M_I=3000$, $M_B=2400$, and $M_\Gamma=816$). In the left panel of Table~\ref{tab:ex3_error_compare}, we present the $L^\infty$ errors of $P_\mathcal{N}$ and $\bU_\mathcal{N}$ obtained using the proposed method under three different $(N_p, N_u)$ pairs. Note that the total number of parameters can be computed as $N_\bft=6N_p+8N_u$. As the network size increases, the proposed method achieves accurate predictions with $L^\infty$ errors ranging from the order of $O(10^{-3})$ to $O(10^{-5})$. The right panel of Table~\ref{tab:ex3_error_compare} lists the results generated from the augmented IIM method~\cite{WT22}. One can see that, the present method can achieve comparably accurate results with the augmented IIM even far less the grid resolutions used. \\

\noindent{\bf Example 4}: As the fourth example (also referred to as Example 6.3 in \cite{WT22}), we keep the domain $\Omega$ and interface $\Gamma$ the same as shown in Example 3. The piecewise constant viscosity is also chosen as $\mu^-=1$ and $\mu^+=0.1$. However, we choose the pressure jump no longer as a constant and the exact solutions of the pressure $p(\bx)$ and velocity  $\bu(\bx)=(u_1(\bx),u_2(\bx),u_3(\bx))$ are given as follows.
\beqs
p(\bx) &=& \left\{ \bary{ll} \left(\frac{3}{8}-\frac{3}{4}x_1^2\right)x_1x_2x_3, & \Omega^-,\\ 0, & \Omega^+,\eary\right.\nonumber\\
u_1(\bx) &=& \left\{ \bary{ll} \frac{1}{4}x_2x_3, & \Omega^-, \\ \frac{1}{4}x_2x_3(x_1^2+x_2^2+x_3^2),  & \Omega^+, \eary\right. \nonumber\\
u_2(\bx) &=& \left\{ \bary{ll} \frac{1}{4}x_3x_1, & \Omega^-, \\ \frac{1}{4}x_3x_1(x_1^2+x_2^2+x_3^2), & \Omega^+, \eary\right. \nonumber\\
u_3(\bx) &=& \left\{ \bary{ll} -\frac{1}{2}x_1x_2(1-x_1^2-x_2^2), & \Omega^-, \\ -\frac{1}{2}x_1x_2x_3^2,  & \Omega^+, \eary\right.\nonumber
\eeqs
and the external force $\bg$ is given
\beqs
g_1(\bx) &=& \left\{ \bary{ll} \left(\frac{3}{8}-\frac{9}{4}x_1^2\right)x_2x_3, & \Omega^-,\\ -\frac{7}{2}\mu^+x_2x_3, & \Omega^+,
	\eary\right. \nonumber \\
g_2(\bx) &=& \left\{ \bary{ll} \left(\frac{3}{8}-\frac{3}{4}x_1^2\right)x_3x_1, & \Omega^-,\\ -\frac{7}{2}\mu^+x_3x_1, & \Omega^+,
	\eary\right. \nonumber \\
g_3(\bx) &=& \left\{ \bary{ll} (\frac{3}{8}-\frac{3}{4}x_1^2-6\mu^-)x_1x_2, & \Omega^-,\\ \mu^+x_1x_2, & \Omega^+. \eary\right.\nonumber
\eeqs
Therefore, the interfacial force has the  form $\bF=(F_1, F_2, F_3)$ as
\beqs
F_1 &=&	-\left(\frac{3}{8}-\frac{3}{4}x_1^2\right)x_1^2x_2x_3 - \frac{3}{4}\dbblk{\mu}x_2x_3 - \mu^+x_1^2x_2x_3 -\frac{1}{2}x_2x_3\left(\mu^+\left(1-x_3^2\right) + \mu^-\left(x_3^2-2x_1^2\right)\right) \nonumber\\
F_2 &=& -\left(\frac{3}{8}-\frac{3}{4}x_1^2\right)x_1x_2^2x_3 - \frac{3}{4}\dbblk{\mu}x_1x_3 - \mu^+x_1x_2^2x_3 -\frac{1}{2}x_1x_3\left(\mu^+\left(1-x_3^2\right) + \mu^-\left(x_3^2-2x_2^2\right)\right) \nonumber\\
F_3 &=& -\left(\frac{3}{8}-\frac{3}{4}x_1^2\right)x_1x_2x_3^2 - \frac{1}{2}\dbblk{\mu}x_1x_2 + 2\mu^+x_1x_2x_3^2 -\mu^-x_1x_2\left(x_3^2-x_1^2-x_2^2\right). \nonumber
\eeqs

Since the domain and the interface are exactly the same as in Example 3, we use the same number of training points $M$. However, we use slightly different neuron pairs $(N_p, N_u)=(30,75), (40,100), (40,120)$ to generate the solution predictions.  Table~\ref{tab:ex4_error_compare} shows the results for different $(N_p, N_u)$ values. One can see that using far less learnable parameters (maximal $N_\bft=1200$ parameters to be learned) and training points ($M=6216$), we are able to produce accurate results comparable to the augmented IIM proposed in~\cite{WT22}.\\
\begin{table}[h!]
\begin{center}
   \begin{tabular}{cccc|cc}
   \hline
     & & Present & &  &  IIM~\cite{WT22} \\
   \hline
   $(N_p,N_u)$ & $E^\infty_p$ & $E^\infty_\bu$ & $\Loss(\bft)$ & Grid no. &  $E^\infty_p$\qquad\qquad$E^\infty_\bu$ \\
   \hline
   $(30,75)$   & $4.44\times10^{-4}$ & $1.24\times10^{-4}$ & $9.26\times10^{-9}$  & $64^3$  &  $6.93\times10^{-4}$\quad$4.07\times10^{-4}$ \\
   $(40,100)$ & $2.32\times10^{-4}$ & $4.63\times10^{-5}$ & $4.61\times10^{-10}$ & $128^3$ &  $1.88\times10^{-4}$\quad$8.67\times10^{-5}$ \\
   $(40,120)$ & $4.14\times10^{-5}$ & $1.56\times10^{-5}$ & $1.42\times10^{-10}$ & $256^3$ &  $3.72\times10^{-5}$\quad$2.06\times10^{-5}$ \\
 	\hline
    \end{tabular}
	\caption{The $L^\infty$ errors of $P_\mathcal{N}$ and $\bU_\mathcal{N}$ for Example 4. Left panel:  the present network results. Here, $\mu^-=1$ and $\mu^+=0.1$, $M=6216$ training data points. Right panel: IIM  results in  \cite{WT22} with grid resolutions $64^3$, $128^3$, and $256^3$.}
	\label{tab:ex4_error_compare}
    \end{center}	
\end{table}

\noindent{\bf Example 5}: Unlike the previous examples with analytic solutions in regular domains, here, we consider the Stokes interface problem with a piecewise-constant viscosity ($\mu^-=10, \mu^+=1$) in the two-dimensional irregular circular domain $\Omega=\{\bx\in\mathbb{R}^2\:|\:\|\bx\|_2\le2\}$ whose solution is not available. The embedded interface $\Gamma$ is parametrized as the polar function $\bx_\Gamma(\alpha)=r(\alpha)(\cos\alpha,\sin\alpha)$, where $r(\alpha)=1-0.2\cos(3\alpha)$ with $\alpha\in[0,2\pi)$, representing a flower shape with three petals. We use $\phi(\bx) = x_1^2+x_2^2-r(\alpha_\bx)$ as the level set function, where $\alpha_\bx=\tan^{-1}\left(\frac{x_2}{x_1}\right)$ so  the augmented feature inputs, $I(\bx)$ and $\phi_a(\bx)$, can be obtained as before accordingly. Along the interface, we now define the surface tension force $\bF(\alpha)=-\gamma\kappa(\alpha)\bn(\alpha)$ to mimic the stationary two-phase flow at some time instant. Here, $\gamma$ is the surface tension and $\kappa$ is the signed curvature of the curve. For simplicity, we set the surface tension $\gamma=1$, the external force $\bg={\bf 0}$, and the boundary condition  $\bu_b={\bf0}$.

\begin{figure}[h]
	\begin{center}
	\includegraphics[width=\textwidth]{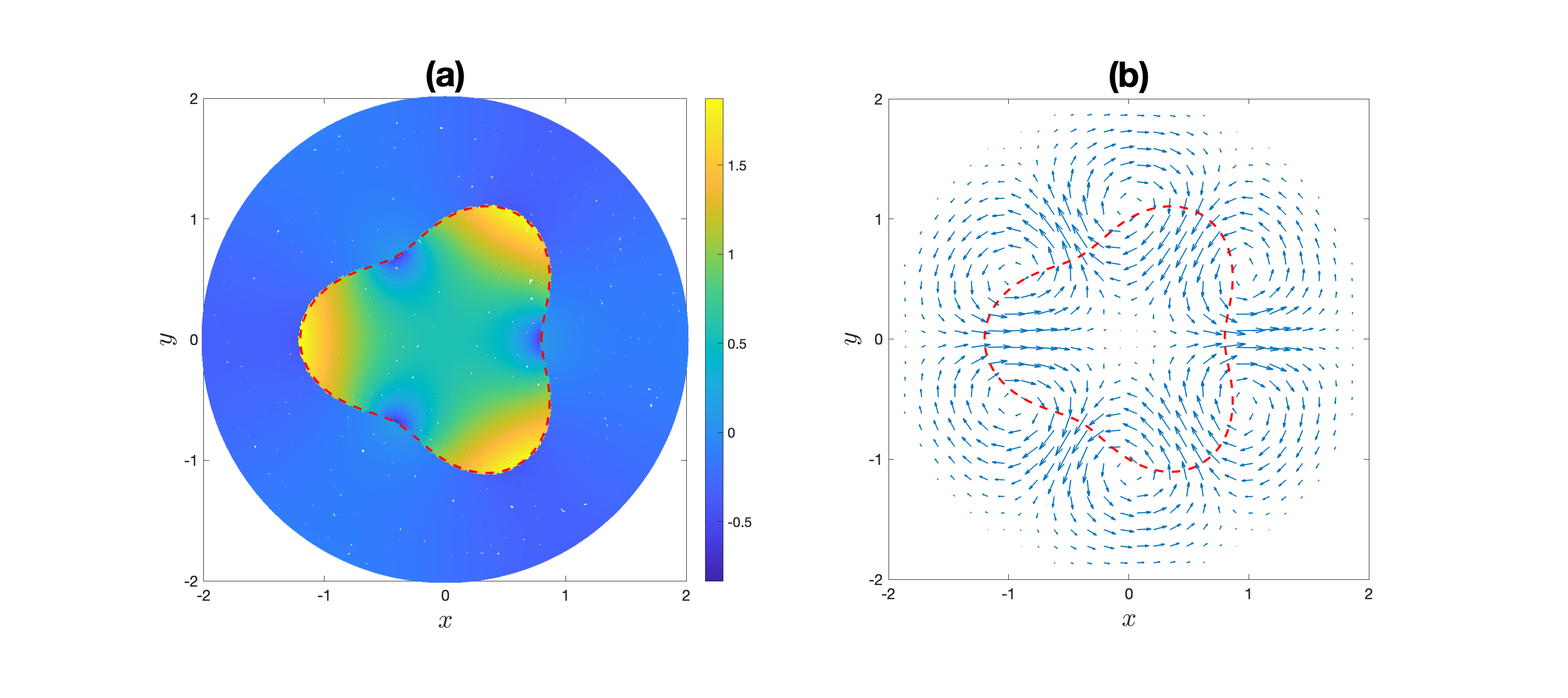}
	\caption{(a) Pressure prediction $P_\mathcal{N}$; (b) Velocity quiver plot of $\bU_\mathcal{N}$.}\label{fig:ex5_mu10}
	\end{center}
	\end{figure}

Since the interface is more complex and the solution is more complicated, we hereby employ a four-hidden-layer deep neural network with $(N_p, N_u) = (15, 20)$ in each hidden layer, resulting in a total number of $N_\bft= 2175$ trainable parameters used. (A comparison with a shallow network will be provided later.) The overall number of training points is chosen as $M=7816$, consisting of $M_I = 4816$ points in $\Omega^-\cup\Omega^+$, $M_\Gamma = 1000$ points on the interface $\Gamma$, and $M_B = 2000$ points on the boundary $\partial\Omega$. We finish the training by either the loss falling below a threshold $\epsilon_\theta=10^{-8}$ or reaching the maximum  $Epoch_{max}=5000$  epochs. Since the analytic solution is unavailable, we investigate the flow behavior by plotting the pressure value and velocity quiver over the whole domain $\Omega$ in Figure~\ref{fig:ex5_mu10}.

Figure~\ref{fig:ex5_mu10}(a) shows the pressure prediction $P_\mathcal{N}$ in $\Omega$ with the interface labeled by the red dashed line. Since the surface tension force acts in the normal direction along the interface, one can clearly see the pressure jump across the interface. The pressure jump is significantly larger, especially near the tip of each petal, where the curvature is relatively large. Meanwhile, near the tip part of each petal (positive curvature region), the pressure inside is higher than the outside one. On the contrary, the pressure inside is lower than the outside one near the neck part (negative curvature region). Figure~\ref{fig:ex5_mu10}(b) displays the corresponding velocity quiver $\bU_\mathcal{N}$ in the domain $\Omega$. It is observed that the flow tends to reduce the absolute magnitudes of curvature along the interface and relaxes to a circular shape. This instant flow tendency matches well with the simulation when the interface dynamics is taken into account. Meanwhile, since the fluid is incompressible, we can see three vortex dipoles (or six counter-rotating vortices) appearing  at the petals.

\begin{figure}[h]
	\begin{center}
	\includegraphics[width=\textwidth]{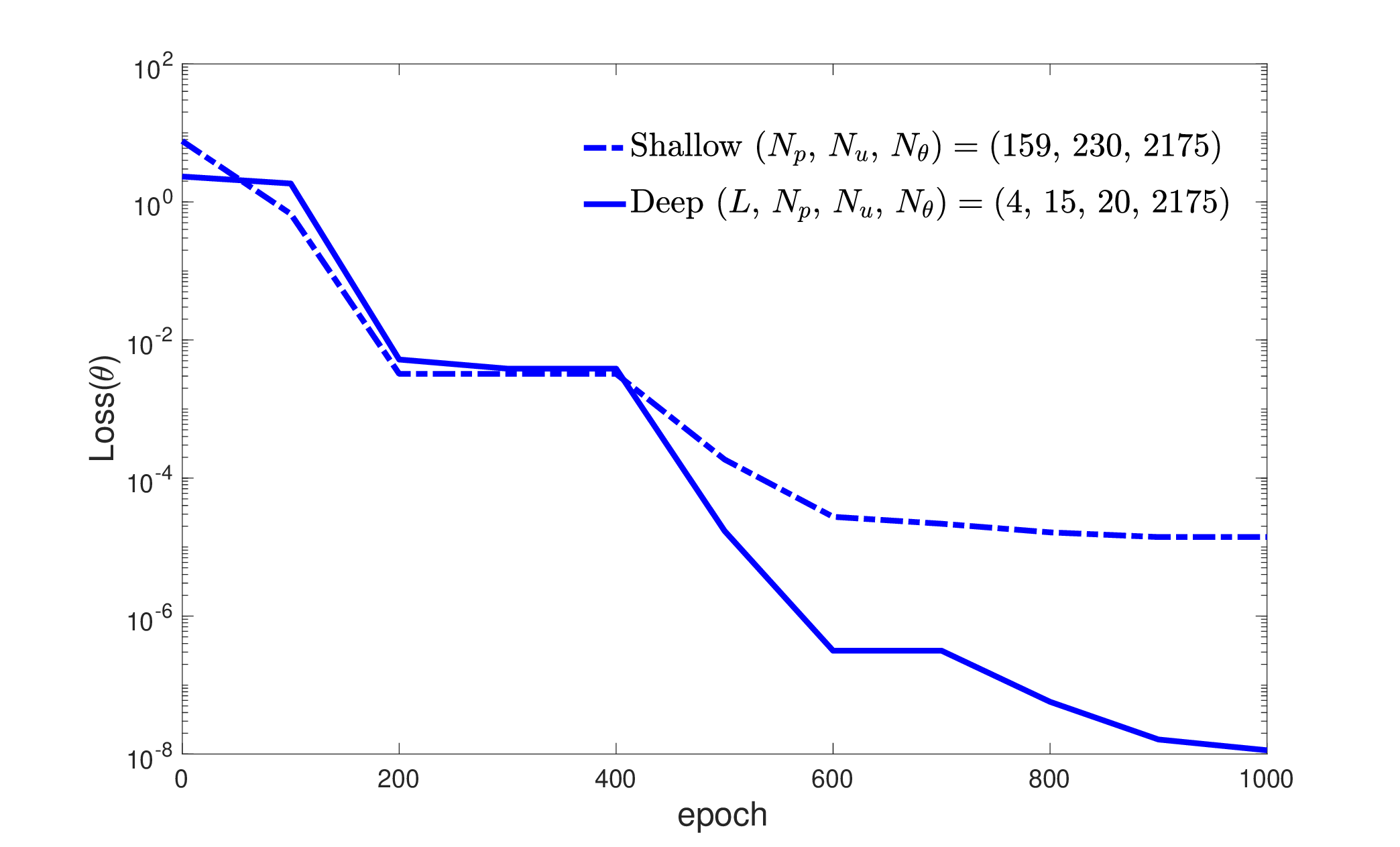}
	\caption{Evolution of training losses using shallow (dash-dotted line) and deep (solid line) networks with the same numbers of trainable parameters $(N_\theta=2175)$ and training points $(M=7816)$.}\label{fig:ex5_loss}
	\end{center}
	\end{figure}

We conclude this example by showing that the deeper neural network used here indeed performs better than a shallow one. To demonstrate that, we use a shallow network structure with neurons $(N_p, N_u) = (159, 230)$ so the total number of trainable parameters $N_\theta = 2175$ is the same as the one used in present four-hidden-layer network with $(N_p, N_u) = (15, 20)$ in each hidden layer. Also we fix the training points $M = 7816$ and train both networks up to $1000$ epochs. The training time for both networks are almost the same since we use the same numbers of trainable parameters and training points. Figure~\ref{fig:ex5_loss} shows both evolution of training losses. One can see that,  using the four-hidden-layer network, the training loss can drop down to the order of $10^{-8}$ within $1000$ epochs, while the training loss of the shallow one only drops to the order of $10^{-5}$. This illustration shows an appropriate network architecture can indeed impact the efficiency of the training process.

\section{Conclusion}
\label{sec:conclusion}
Solving Stokes equations with a piecewise-constant viscosity and singular force on an interface using traditional grid-based methods encounter the major challenge of computing the non-smooth solutions numerically. One can reformulate the governing equations by solving the Stokes equations in two sub-domains separately; however, the pressure and velocity from both sub-domains are coupled together by the traction balance equation on the interface. This makes the grid-based methods hard to discretize accurately near the interface and the resultant linear system hard to solve. In this paper, we provide a neural network approach to solve the equations. Our approach is based on adopting a discontinuity capturing sub-network to approximate the pressure and a cusp-capturing sub-network to approximate the velocity so that the network predictions can retain the inherent properties of the solutions. Since the present neural network method is completely mesh-free, the implementation is much easier than the traditional grid-based methods, such as immersed interface method, which needs careful and laborious efforts to discretize the equations accurately near the interface. A series of numerical experiments are performed to test the accuracy against the results obtained from the immersed interface methods in the existing literature. A completely shallow (one-hidden-layer) network with a moderate number of neurons and sufficient training data points can obtain comparably accurate results obtained from the traditional methods. We shall apply the present network to solve realistic two-phase flow applications in the future.

\section*{Acknowledgments}
Y.-H. Tseng and M.-C. Lai acknowledge the supports by National Science and Technology Council, Taiwan, under research grants 111-2115-M-390-002 and 110-2115-M-A49-011-MY3, respectively.


\end{document}